\documentclass[12pt]{article}
\usepackage{amsmath,amsthm,amsfonts,amssymb}
\usepackage{verbatim}
\usepackage{latexsym}
\newtheorem{theorem}{Theorem}[section]
\newtheorem{lemma}[theorem]{Lemma}

\newtheorem{definition}[theorem]{Definition}

\newtheorem{remark}[theorem]{Remark}

\newcommand\NN{{\Bbb N}}
\newcommand\Z{{\Bbb Z}}
\newcommand\RR{{\Bbb R}}

\newcommand\ZZ{{\Bbb Z}}

\begin{document}
\title{ Asymptotic Uncorrelation for Mexican Needlets}
 \author{
  \thanks{ This work was partially  supported by the Marie Curie Excellence Team Grant MEXT-CT-2004-013477,
   Acronym MAMEBIA, funded by the European Commission.}
Azita Mayeli \\
\footnotesize\texttt{{amayeli@math.sunysb.edu}} }

\date{\today}
\maketitle

\begin{abstract}
We recall Mexican needlets from \cite{gm2}. We derive an estimate
for certain types of Legendre series, which we apply to the
statistical properties of Mexican needlets. More precisely, we shall
show that, under isotropy and Gaussianity assumptions, the Mexican
needlet coefficients of a random field on the sphere are
asymptotically uncorrelated, as the frequency parameter goes to
infinity.  This property is important in the analysis of cosmic
microwave background radiation.
\end{abstract}

\vspace{.5cm}

\footnotesize{
\begin{tabular}{lrl}
{\bf  Keywords.} & {\em Spherical harmonics, spherical Laplacian operator,  wavelets, Mexican  needlets,}  \\
   &  \multicolumn{2}{l} {\em   angular power spectrum, Legendre (Gegenbauer) polynomials.
  }\\
  %&  \hspace{-.1cm}  {\em  }\\
 \end{tabular}
 \vspace{.3cm}

 \begin{tabular}{lrl}
  &&  \hspace{-1cm}{\bf  AMS  Subject Classification (2000).} {62M40, 33C55, 65T60}
\end{tabular}
}

%INTRODUCTION
\section{Introduction}\label{introduction}

 Let  $\Delta$  denote  the spherical Laplacian. Let
  $\{Y_l^m\}$ for $l\in \NN$ and $m= -l, \cdots, l$  be the usual spherical harmonics on $S^2$.  They constitute an
orthonormal basis for $L^2(S^2)$ and are eigenfunctions for the
Laplacian operator with eigenvalues $l(l+1)$, i.e.,
 $\Delta Y_l^m = l(l+1)Y_l^m$. We shall denote $\lambda_l= l(l+1)$.\\

By  the spectral theorem for  the spherical Laplacian $\Delta$, for
any  $f\in \mathcal{S}(\mathcal{\RR^+})$ and $t>0$, the operator
 $f(t^2\Delta)$ is a bounded operator on $L^2(S^2)$. Let
 $K_t$ denote the kernel of operator $f(t^2\Delta)$ defined on $S^2\times S^2$. In \cite{gm2} we
observed that $f(t^2\Delta)$ has an explicit kernel $K_t$, given by
\begin{align}\label{kernel}
 K_t(x,y) &= \sum_{l=0}^\infty \sum_{m=-l}^l f(t^2\lambda_l)  Y_l^m(x) \overline{ Y_l^m}(y).
 \end{align}
It is easy to see that $K_t$ is smooth in $(t,x,y)$ for $t > 0$.

In applications to wavelets, one needs to assume $f(0)=0$, and we do
so in what follows.

The localization properties of these kinds of kernels have already
been studied in detail in \cite{gm2}; briefly,  we showed that for
every pair of $C^\infty$ differential operators $X$ (in x) and $Y$
(in y) on $S^2$, and for every integer $N\in \NN_0$, there
 exists $c:=c_{N,X,Y}$ such that for all $t>0$ and $x, y\in S^2$
\begin{align}\label{gaussian-decay}
\left| XYK_t(x,y) \right|\leq c ~ \frac{ t^{-(2+I+J)}}{\left(\frac{
d(x, y)}{t}\right)^N},
\end{align}
 where $I:=\deg X$ and $J:=\deg Y$.  (In the case where $f$ has compact support away from $0$,
and if one replaces $\lambda_l=l(l+1)$ by $l^2$ in the formula for
$K_t$, this was earlier shown by Narcowich, Petrushev and Ward, in
\cite{Narcowich06} and \cite{Narcowich2}.)  In fact our argument
worked on general smooth compact oriented Riemannian manifolds, not
just the sphere.

Moreover, on the sphere $S^2$, the kernel  $K_t(x,y)$ in
(\ref{kernel}) is rotationally invariant, i.e.,
 $K_t(x,y)= K_t(\rho x, \rho y)$
for any rotation $\rho$ defined on the sphere. Intuitively, we
define  $\psi_{t,x}(y):= K_t(x,y)$, then  $\psi_{t,x}(y)$  can be
thought of as analogous to ``$t$-dilation" and  ``$x$-translation"
of the single function $\psi(y):= K_1({\bf 1},y)$, where ${\bf
1}=(1,0,0)$ is the ``north pole". For more about the analogues of
``dilation" and ``translation"  we refer the reader to our earlier
work
\cite{gm2}.\\

 In \cite{gm2} we  also proved  that
 if  $f\in \mathcal{S}(\mathcal{\RR^+})$ with $f(0)=0$
  satisfies a discrete version of Calder\'on's formula, then
  the associated kernel is a wavelet, i.e.,   for $a>1$ and
 for a carefully   chosen discrete set  $\{x_{j,k}\}_{(j,k)\in \Z\times \Z}$ on the sphere
and certain weights $\mu_{j,k}$, the collection of $\{\psi_{j,k}:=
\mu_{j,k}K_{a^j}(x_{j,k}, \cdot)\}$
 constitutes a wavelet  frame.  By this we  mean that    there exist constants $0<A\leq B<\infty$,
such that  for any $F\in L^2(S^2)$ the following holds:
 \begin{align}\label{frame}
 A\parallel F\parallel_2^2 \leq \sum_{j,k} \mid \langle F, \psi_{j,k}\rangle\mid^2 \leq B\parallel F\parallel_2^2.
 \end{align}
 We also showed that, if the points $\{x_{j,k}\}$
were selected carefully enough, $B$ could be made arbitrarily close
to $B_a$ and $A$ could be made arbitrarily close to $A_a$, where
 $B_a/A_a \rightarrow 1$ almost quadratically as  $a\rightarrow 1$.
The frame is therefore a ``nearly tight frame". Like before,
$\psi_{j,k}$ is analogous to $a^j$- dilation and
$x_{j,k}$-translation of the function $K({\bf 1}, \cdot)$.  We shall
call such a frame needlets based on $f$.  (The term ``needlets'' was
used earlier for the frames of Narcowich, Petrushev and Ward.)

 %MAIN CONTRIBUTIONS
%\subsection{Main Contributions}

The present work involves the  needlets  based on the Schwartz
functions
\begin{equation}
\label{srf0} f(s)= s^r f_0(s),
\end{equation}
where  $0 \neq f_0\in \mathcal{S}(\RR^+)$ and $r\in \NN$. For
$f_0(s)= e^{-s}$ and $r\in \NN$,
we shall call these needlets  {\bf Mexican needlets}.\\

Let $ \psi^{t,x}= K_t(x,\cdot)$, so that
$\psi_{j,k}=\mu_{j,k}\psi^{a^j, x_{j,k}}$.
 By the spherical harmonic  expansion of the kernel of $f(t^2\Delta)$ in (\ref{kernel}),   for
any $F\in L^2(S^2)$ at every scale $t$ and every  pixel  $x$   we
have
\begin{align}\label{ran-need-coef}
\beta_{t, x}:= \langle F,  \psi^{t,x} \rangle= \sum_{l=1}^\infty
\sum_{m=-l}^l f(t^2\lambda_l) \langle F, Y_l^m\rangle  ~Y_l^m(x).
\end{align}
It is standard to  denote $\hat{F}(l,m)=\langle F,Y_{l,m}\rangle$
and call it the spherical harmonic coefficient
 of the function $F$ at $(l,m)$.

In applications to cosmology (\cite{BKMP06-3}), one assumes that $F$
is a centered  Gaussian field which is mean square continuous and
isotropic (which means that $E(F(\rho(x)\rho(y)) = E(F(x)F(y))$ for
every orthogonal transformation $\rho$). In that case, if we set
$a_{l,m}=\hat{F}(l,m)$, then the $a_{l,m}$ are a set of  Gaussian
random variables, uncorrelated except for the condition that
$\overline{a}_{lm} = (-1)^m a_{l,-m}$.  We define

\begin{align}
\beta_{t, x}= \langle F, \psi^{t,x}\rangle=\sum_{l=1}^\infty
\sum_{m=-l}^l f(t^2\lambda_l) a_{l,m}  ~Y_l^m(x).
\end{align}
Note $\beta_{t, x}$ is real, since $\overline{Y_{l,m}}= (-1)^m
Y_{l,-m}$. Taking $t=a^j$ and $x=x_{j,k}$, ~$\mu_{j,k}\beta_{t,x}$
become the random needlet coefficients $\langle F,
\psi_{j,k}\rangle$. \\

The isotropy assumption implies that $E(a_{l,m}a_{l',m'})=0$ for
$l\neq l', ~m\neq m'$ and $E(a_{l,m}^2)=c_l$ (independent of $m$).
Then for any $t>0$ and  $x,y\in S^2$,
   the expectation of $ \beta_{t, x}\beta_{t, y}$ is   given by
 \begin{align}\label{expectation-x-y}
 E( \beta_{t, x}\beta_{t, y})= \sum_l f(t^2\lambda_l)^2 c_l (2l+1) P_l^{1/2}(x\cdot y),
 \end{align}
where $P_l^{1/2}$ is the Legendre polynomial of degree $l$ (or
 Legendre polynomial of index $1/2$ and degree $l$) and $x\cdot y$ is
the usual inner product of $x$ and $y$. The $P_l^{1/2}$ may be
defined in terms of the generating function
               \begin{align}
               (1-2\xi\eta+\xi^2)^{-1/2}= \sum_{l=0}^\infty P_l^{1/2}(\eta)~\xi^l.
               \end{align}

            The main goal of this work is to study the behavior of
the following quantity, for different $x, y$ on the sphere, when $t$
is close to zero:
            \begin{align}\label{rational-term}
     \text{Cor}(\beta_{t,x}, \beta_{t,y})=  \frac{E( \beta_{t, x}\beta_{t, y})}{ \sqrt{E( \beta_{t, x}^2 )} \sqrt{E( \beta_{t, y}^2 )}}.
                      \end{align}

From the statistical point of view, (\ref{rational-term}) is
refereed to the correlation of the
           $ \beta_{t, x}$ and $\beta_{t, y}$ at scale $t$ and at positions  $x, y$ respectively and is denoted by
  $\text{Cor}(\beta_{t, x},\beta_{t, y})$.
           Using  this statistical terminology, we shall show that the
           correlation
approaches zero as $t \rightarrow 0$, if $x$ and $y$ are fixed.
That is, they are asymptotically  uncorrelated. For this, first  we
prove  in
          Lemma \ref{tool} that
             the term in  (\ref{expectation-x-y}) is well localized near $x=y$.
 Then applying this lemma, we  present the results  concerning the
asymptotic  uncorrelation  in   Theorem  \ref{main-theorem}. Note
that if $p,q>0$, then $ \text{Cor}(p\beta_{t, x},q\beta_{t, y})=
 \text{Cor}(\beta_{t, x},\beta_{t, y})$. Thus if $t=a^j$,
 $x=x_{j,k}$, $y=x_{j,k'}$, for certain $k,k'$, the correlation in
 (\ref{rational-term})  gives the correlation between needlet
 coefficients. \\

Although the motivation is from statistics and cosmology, our
arguments to prove Theorem  \ref{main-theorem} will be purely
mathematical and will not use this motivation. However, in order to
prove this theorem, we shall assume some regularity assumptions for
the expected values $c_l$. (These values are known as the angular
power spectrum in astrophysics.) Our assumptions will be reasonable,
based on the astrophysical literature.
   More precisely,
 we  suppose that $c_l$ is given by  formula of a  function $u\in C^\infty(\RR^+)$, $c_l:=u(l)$,
 with the following properties:
For some  real number $\alpha>2$:\\
(a) for any $k\in \NN_0$ there exists a constant $c_k$ such that

\begin{align}\label{alpha}
|\partial^k u(s)|\leq C_k~ s^{\alpha-k} \quad \forall ~ s\geq 1,
\end{align}
and   (b)  there exist positive constants $c_0, c_1$  such that for
any $s\in (1,\infty)$,
$c_0 s^{\alpha} \leq u(s)\leq c_1 s^{\alpha}$ holds.  \\

 With these  assumed conditions for  the angular power spectrums
 $c_l$, in Theorem \ref{main-theorem}
                        we will show  that
         for any  fixed points $x$ and $y$ the expression for the correlation
         of
   $ \beta_{t, x}$ and $\beta_{t, y}$  satisfies
   \footnote[1]{
While finishing this paper, we learned by personal communication
that simultaneously   X. Lan and D. Marinucci \cite{lan-marri} have
obtained an analogous result to Theorem \ref{main-theorem}. The
  assumptions are not
equivalent and the approaches are entirely different. We believe
both are of independent interest and should be utilized.}: there
exists a constant $C > 0$ such that
               \begin{align}\label{good-loc}
               |\text{Cor}(\beta_{t, x},\beta_{t, y})|=
 \frac{|E( \beta_{t, x}\beta_{t, y})|}{ \sqrt{E( \beta_{t, x}^2 )} \sqrt{E( \beta_{t, y}^2 )}}\leq
\frac{C t^{4r-\alpha+2}}{\left({d(x,y)}\right)^{2N}},
               \end{align}
uniformly in $(t,x,y)$, provided $4r-\alpha+2 > 0$.  (Here, $r$ is
as in (\ref{srf0}) and $N$ is the least integer greater than
$2r-\alpha/2 + 1$.) Hence, the asymptotic uncorrelation of
coefficients for small $t$ holds, for $r$ sufficiently large in
relation to $\alpha$.
 We shall prove this in   Theorem \ref{main-theorem}.\\

Our main contributions in this work will be as follows: After this
introduction to our work,
  which will be concluded by some further historical comments about needlets and Mexican needlets,
in Section \ref{Not-Pre} we will present some notation and
preliminary definitions.
 In Section \ref{motivation} we shall present motivation for  this work.
In section \ref{the-key-estimate}, first we study   the localization
property of  (\ref{expectation-x-y})
 in  Lemma \ref{tool},  which has a crucial role in the proof of the main theorem.
 Applying this lemma,   in Theorem   \ref{main-theorem} we prove that the random needlet coefficients,
 with the assumed properties for the $c_l$,  are asymptotically uncorrelated  for any fixed  angle  on the sphere,
for $r$ sufficiently large in relation to $\alpha$.

  \subsection{More about Needlet and Mexican Needlets: History}\label{history}

Needlets were introduced by Narcowich, Petrushev and Ward in
\cite{Narcowich06}. In place of our function $f$, they consider only
smooth functions $g$ with compact support away from $0$.  We prefer
to call such a $g$ a ``cutoff function".  As we have indicated, they
also used $l^2$ in place of $\lambda_l = l(l+1)$ in the definition
(see (\ref{kernel})), which is a minor distinction.

Asymptotic uncorrelation was first studied in \cite{BKMP06-3}.
There the authors assumed that $c_l:= \frac{g_j(a^{-j}l)}{l^\alpha}$
for every $l$ such that $a^j < l < a^{j+1}$; here
 $(g_j)_j$ is a sequence of functions (they can also be constants), which have a uniform bounded
differentiability condition up  to order $M$, for some large $M$.
The investigation of needlets from a stochastic point of view is
due to \cite{BKMP07, BKMP06-2, BKMP06-3}. Needlets  have  been used
by astrophysicists to study cosmic microwave background radiation
(CMB).
(See for instance \cite{baldi, guil} and references therein.)  \\

Needlets and Mexican needlets each have their own advantages.
Needlets have these advantages: for appropriate $f$ and $x_{j,k}$,
needlets are a tight frame on the sphere (i.e., $A=B$ in
(\ref{frame})). The frame elements at non-adjacent scales are
orthogonal.  The random needlet coefficients are asymptotically
uncorrelated  (there is no ``$r$'' that needs to be assumed large in
relation to $\alpha$).

Mexican needlets (as developed in \cite{gm2}, for which $f(s) = s^r
e^{-s}$) have their own advantages. We write down an approximate
formula for them which can be used directly on the sphere in
\cite{gm2}. (This formula, which arises from computation of a
Maclaurin series, has been checked numerically. It is work in
progress, expected to be completed soon, to estimate the remainder
terms in this Maclaurin series.)  Assuming this formula, Mexican
needlets have Gaussian decay at each scale. They do not oscillate
(for small $r$), so they can be implemented directly on the sphere,
which is desirable if there is missing data (such as the ``sky cut''
of the CMB).  Finally, as the proofs in this article will show, the
constant $C$ in (\ref{good-loc}) depends on finitely many
derivatives of $f$; so in order to be sure that this constant is
small as possible, it is desirable to use real analytic functions,
such as $s^r e^{-s}$, as $f(s)$.

In our opinion, both needlets and Mexican needlets should be
utilized in the analysis of CMB, and the results should be compared.

%NOTATIONS AND PRELIMINARIES
\section{Notations and Preliminaries}\label{Not-Pre}

Let $S^2$ denote the unit sphere in $\RR^3$ and let ${\bf 1}= (1, 0,
0)$ denote the  ``north pole". We shall consider $S^2$ with its
rotationally invariant measure $\mu$. We may write
\begin{align}
L^2(S^2)=\oplus_{l=0}^\infty \mathcal{H}_l,
\end{align}
where $\mathcal{H}_l$ is the space of spherical harmonic of degree
$l$. In fact, $\mathcal{H}_l$ is the restriction of homogeneous
harmonic polynomials of degree $l$ on $\RR^3$. The Laplacian
operator for the sphere is in a sense the ``restriction" of the
Laplacian operator for $\RR^3$ on $S^2$. We shall denote it by
$\Delta$. A precise formula for $\Delta$  in spherical coordinates
is given by
\begin{align}
\Delta= \frac{1}{\sin\theta}\frac{\partial}{\partial\theta}\left(
\sin\theta\frac{\partial}{\partial\theta}\right) +
\frac{1}{\sin^2\theta}\frac{\partial^2}{\partial\varphi^2}
\end{align}
where  $0 \leq  \theta \leq \pi $ and
$0 \leq \varphi < 2\pi$.  \\
   If $P \in
{\mathcal H}_l$, then
\[ \Delta P = l(l+1) P,\]
which means that the spherical harmonics of degree $l$ are eigenfunctions for $\Delta$ with eigenvalues $\lambda_l= l(l+1)$.\\

Within each space ${\mathcal H}_l$ is a unique {\em zonal harmonic}
$Z_l$, which has the property that for all $P \in {\mathcal H}_l$,
$P({\bf 1}) =  \langle  P , Z_l  \rangle $. In particular, $P$ is
orthogonal to $Z_l$ if and only if $P({\bf 1}) = 0$. Obviously,
$Z_l({\bf 1}) =  \langle  Z_l,Z_l  \rangle$ and  $Z_l(y)$ is known
explicitly in terms of the Legendre polynomials.  In fact, if
$\omega_2$ is the area of $S^2$, then for $c= \frac{1}{\omega_2}$,

\begin{equation}
\label{zony} Z_l(y) = c (2l+1) P_l^{\lambda}(y_1),
\end{equation}
  where $y = (y_1,y_2, y_3)$, $\lambda = 1/2$, and
$P^{\lambda}_l$ is the Legendre polynomial of degree $l$ associated
with $\lambda$ \cite{Stein-Weiss71}. In the sequel we shall avoid
the notation $\lambda$ as $\lambda=1/2$. Since
\begin{equation}
\label{pk1} P_l (1) = \left(\displaystyle^l_l\right) = 1,
\end{equation}
then
\begin{equation}
\label{zony-one} Z_l({\bf 1}) = c  (2l+1).
\end{equation}

$\mathcal{H}_l$ is a finite dimensional vector space in $L^2(S^2)$
with dimension
$$\dim \mathcal{H}_l= 2l+1.$$
We shall choose an orthonormal basis for each ${\mathcal H}_l$, one
of whose elements is $Z_l/\|Z_l\|_2$.  For $l\in \NN_0$, let
$\{Y_l^m\}_{m=-l}^l$ be an orthonormal basis of $2l+1$ elements for
$\mathcal{H}_l$ with central element $Y_l^0= Z_l/\|Z_l\|_2$, and
with $\overline{Y_m^l} = Y_{-m}^l$. Therefore  any $F\in L^2(S^2)$
has a spherical harmonic expansion  $\{Y_l^m\}_{m=-l}^l$:
\begin{align}\label{spherical-expansion}
F= \sum_{l=0}^\infty \sum_{m=-l}^l  \hat F(l,m) ~Y_l^m
\end{align}
where $\hat F(l,m):=
\langle F, Y_l^m\rangle$ are called the spherical harmonic coefficients. \\

Next we introduce notation concerning the discrete version of
derivation on a sequence:

 \begin{definition} Let  $\Delta^+$ and $\Delta^-$  denote
the   difference operators   defined on
  any sequence
  $\{a_l\}_{l\in \ZZ}$ as follows: For all $l$,
  \begin{align}
\Delta^+ a_l =& a_{l+1}-a_l\\
\Delta^- a_l=& a_l-a_{l-1}.
\end{align}
Obviously, the operators  $\Delta^-$ and $\Delta^+$ commute.
 \end{definition}
Suppose $\{b_l\}_{l\in \ZZ}$ and $\{c_l\}_{l\in \ZZ}$ are any pair
of sequences. Then the product rules for the difference operators
are given as follows:
\begin{align}
\Delta^+(b_l c_l)&= (\Delta^+b_l)c_{l+1}+b_l(\Delta^+c_l),\\
\Delta^-(b_lc_l)&= (\Delta^- b_l ) c_{l}+b_{l-1} (\Delta^-c_l).
\end{align}

For two sequences   $\{a_l\}_{l\in \ZZ}$   and   $\{b_l\}_{l\in
\ZZ}$, we shall say $a_l$ is uniformly bounded by $b_l$ from above,
and write $a_l=\mathcal{O}(b_l)$, if for some positive constant $c$
we have
 $$| a_l |\leq c~ |b_l|\quad \forall ~l\in \ZZ.$$

% MOTIVATION
\section{Motivation}\label{motivation}

Following the notations  of our earlier work in  \cite{gm2}, for any
$f\in \mathcal{S}(\RR^+)$ with $f(0)=0$ and $t>0$, let $K_t$ be the
associated kernel of the operator $f(t^2\Delta)$, i.e., for any
$F\in L^2(S^2)$
\begin{align}
f(t^2\Delta)F(x)= \int_{S^2} K_t(x,y)F(y)~d\mu.
\end{align}

By expansion in   (\ref{spherical-expansion}) and the spectral
theorem,  one has:
\begin{align}
\label{kersph1}
\psi_{t, x}(y):=K_t(x,y) &= \sum_{l=1}^\infty \sum_{m=-l}^l f(t^2\lambda_l)  Y_l^m(x) \overline{ Y_l^m}(y)\\
&= \frac{1}{\omega_2}~\sum_{l=1}^\infty  f(t^2\lambda_l) ~ (2l+1)
P_l(x.y).
 \end{align}
 From now on, we shall neglect the factor $\frac{1}{\omega_2}$ in the rest of the  work, since it will   not affect  our main results. \\

 Since $\langle Y_l^m, Z_l\rangle= Y_l^m({\bf 1})= 0$, then
 \begin{align}
\psi_t(y):=\psi^{t,{\bf 1}}(y) = \sum_{l=1}^\infty f(t^2\lambda_l)
\overline{ Z_l}(y)=  \sum_{l=1}^\infty f(t^2\lambda_l)   (2l+1)
P_l(y_1).
 \end{align}

 As we explained in the introduction, the
   $\psi_t$, generates a needlet frame (based on $f$) for $L^2(S^2)$.\\

 For  the scale $t>0$ and position $x\in S^2$,
 \begin{align}
\beta_{t, x}=\sum_{l=1}^\infty \sum_{m=-l}^l f(t^2\lambda_l)
\langle F, Y_l^m\rangle  Y_l^m(x).
 \end{align}

% \begin{definition}[Mexican Needlets] ~
 %For $f_0\in \mathcal{S}(\RR^+)$,  $m\in \NN$,  let $f(s)= s^r f_0(s)$, the spherical   wavelets  defined    in %(\ref{kersph1})
%are    called {\bf  mexican needlets}, and, for $x, t$,  the coefficient  $\beta_{t,x}$ is called the mexican %needlte coefficient. (The name of mexican needlets is  derived  from the name of
 % the mexican hat wavelet   introduced  in \cite{gm2},  for which
    %$m=1$ and $f_0(s)=e^{-s}$.)
 %\end{definition}

 Now, suppose that $\{a_{l,m}\}$ are
 Gaussian random variables.   For $t>0$ and $x\in S^2$,  we set

  \begin{align}
\beta_{t, x}=\sum_{l=1}^\infty \sum_{m=-l}^l f(t^2\lambda_l) a_{l,m}
Y_l^m(x)
 \end{align}
For $t>0$ and any other point  $y\in S^2$,  we define the
correlation:
  \begin{align}\label{un-corr}
 \text{Cor}(\beta_{t,x}, \beta_{t,y})=\frac{E(\beta_{t,x}\beta_{t,y})}{
 \sqrt{E(\beta_{t,x}^2)} \sqrt{E(\beta_{t,y}^2)}
 }
 \end{align}
  where for any random variable $w$,  $E(w)$ is its  expectation.\\

 As motivated in the introduction, we assume that there exist $\{c_l\}$ such that
 \begin{align}\label{expect}
 E(a_{l,m}a_{l',m'})=c_l \delta_{l,l'}\delta_{m,m'}\quad \forall ~~ l, l', ~m, m'\in \ZZ.
 \end{align}

We first  calculate $E(\beta_{t,x}\beta_{t,y})$ for any $t>0$ and $x, y\in S^2$:\\
    By the linearity of $E$ and  using (\ref{expect})  one obtains:
 \begin{align}
 E(\beta_{t,x}\beta_{t,y})%&= E\left(\sum_{l,l', m,m'} a_{l,m}a_{l',m'} f(t^2\lambda_l)f(t^2\lambda_{l'})Y_{l}^m (x)
%Y_{l'}^{m'}(y)\right)\\
 &=\sum_{l, m}  E(a_{l,m}^2) f(t^2\lambda_l)^2 Y_{l}^m (x)
 Y_{l}^{m}(y)
 \end{align}
 Surely  $c_l:=E(a_{l,m}^2)$ (which is called  the angular power spectrum).  Therefore
  \begin{align}
 E(\beta_{t,x}\beta_{t,y})= \sum_l (2l+1) c_l f(t^2\lambda_l)^2  P_l(x.y),
 \end{align}
 and, hence,    for  $x=y$ we obtain

  \begin{align}
 E(\beta_{t,x}^2)= \sum_l (2l+1)   c_l f(t^2\lambda_l)^2  P_l(1).
 \end{align}
 Therefore the correlation formula   is  given by
 \begin{align}\label{un-corr-2}
 \text{Cor}(\beta_{t,x}, \beta_{t,y})=\frac{ \sum_l  (2l+1)  c_l  f(t^2\lambda_l)^2  P_l(x.y)
 }{
 \sum_l   (2l+1)  c_l  f(t^2\lambda_l)^2  P_l(1)
 }.
 \end{align}

A similar formula was derived in \cite{BKMP06-3}.

We will
to estimate (\ref{un-corr-2})  for small  $t> 0$  with  following  assumptions on the angular power spectrums $c_l$:\\
Assume that   $c_{l}$  is given by
   \begin{align}
   c_l= u(l)
    \end{align}
    where $u$ is a smooth  function on $(0,\infty)$   which satisfies the  following conditions, for some  $\alpha\in \RR_{> 2}$:
   \begin{itemize}
     \item{(i)} for all $k\in \NN_0$ there exists a constant $C_k$ such that
  $|\partial^ku(s)| \leq C_k~ s^{-\alpha-k}$ (uniformly for all $s\geq 1$), and
     \item{(ii)} there exist $k_0, k_1>0$ such that
       \begin{align} k_0 s^{-\alpha}\leq u(s)\leq k_1 s^{-\alpha }\quad \forall ~ s\geq 1.
       \end{align}
    \end{itemize}

With above assumptions on $c_l$, in the next section we   show that
if  $4r+2> \alpha$, then for  two different points $x, y$ on the
sphere
 and
  $t>0$  one has:

 \begin{align}\label{cor}
\left| \text{Cor}(\beta_{t,x}\beta_{t,y})\right|=\left| \frac{
\sum_{l=1} f(t^2\lambda_l)^2 c_l Z_l(x\cdot y)}{
\sum_{l=1}f(t^2\lambda_l)^2 c_l Z_l({\bf1})}\right| \leq
\frac{C t^{4r-\alpha+2}}{d(x, y)^{2N}},
 \end{align}
 where $C:=C(r,\alpha)$ is independent of choice of  $x, y$ and $t$ and $N$ is the least integer greater
than $2r-\alpha/2+1$. Here $r$ is as in (\ref{srf0}).  The proof
will show that in (i) it is enough to assume the
condition for all $k \leq M_0$ for some large $M_0$.\\

% Note that if $f$ vanishes near zero, we may take $r$ arbitrary, so the correlation
%is less than or equal to $C_M/[\frac{d(x,y)}{t}]^M$ for any $M$.
%This was more or less the situation studied in \cite{BKMP06-3}.
%(They used $l^2$ instead of $\lambda_l = l(l+1)$; this is a
%minor difference.)\\

    Here we shall present an example for $u$ for which  the  conditions (i)-(ii) are satisfied: \\
\ \\
 {\bf  Example:}
Say $u(s) = \frac{F(\log s) P(s)}{s^{\beta}Q(s)}$ on $(0,\infty)$,
where: $P$ and $Q$ are polynomials of degree $p$ and $q$
respectively, and $\beta+q-p = \alpha$; $F$ is a smooth function
defined on $\RR$ with all bounded derivatives; and $P$, $Q$ and $F$
are all positive on $[1,\infty)$, and $F > \tau > 0$ on
$[0,\infty)$.  Then, surely, $u$ satisfies (i) and (ii).

%THE KEY ESTIMATE
\section{The Key Estimate}\label{the-key-estimate}
 For the study
    of the behavior    of   \begin{align}
 \frac{ \sum_l  (2l+1)  c_l  f(t^2\lambda_l)^2  P_l(x.y)
 }{
 \sum_l   (2l+1)  c_l  f(t^2\lambda_l)^2  P_l(1)
 }
 \end{align}
 when $t$ is close to zero,
we   first prove the next lemma, which has a crucial role in the
proof of our main theorem, Theorem \ref{main-theorem}:

  \begin{lemma}\label{tool}
  Suppose that for
 $\{a_l\}_{\ZZ}$   and $\mu\in \RR$  the  following hold:
 \begin{itemize}
 \item[(i)]   $ a_l= \mathcal{O}(l^{\mu}) \quad \forall l\in \NN$, $a_0=0$, and
  \item[(ii)]   for all $k_1, k_2\in \NN_0$, $  ({\Delta^-})^{k_1}({\Delta^+})^{k_2} a_l =
 \mathcal{O}(l^{\mu-2(k_1+k_2)})$.
   \end{itemize}
   Define
  $f(\cos\theta):= \sum_{l=1}^\infty  a_l~ Z_l(\cos\theta)$. If $\mu+2>0$, then for some positive constant
the following inequality holds uniformly:
\begin{align}\label{main}
\left| f(\cos\theta)\right|\leq \frac{C}{\mid
\theta\mid^{2N}}\quad ~~ \forall ~~0<\theta\leq \pi,
\end{align}
when $N$ is the least integer greater than $\mu/2+1$.  (Observe that
here $Z_l(\cos\theta)$ is understood as \\
$(2l+1)P_l(\cos\theta)$
for the Legendre polynomial of degree $l$.)
\end{lemma}

  \begin{proof}

\begin{comment}
For any $\theta$, write:
\begin{align}\label{spreading}
f(\cos\theta)= \sum_{l\leq b} a_l~ Z_l(\cos\theta)+\sum_{l> b} a_l~
Z_l(\cos\theta).
\end{align}
 For the first summation,
using the  property  (i) for the  coefficients $a_l$ and that
$-1\leq P_l(x)\leq 1$, and since $\mu+2>0$,  we have:
 \begin{align}\label{first-part}
\mid  \sum_{l\leq b}a_l~ Z_l(\cos\theta)\mid \leq \sum_{l\leq
b}(2l+1) l^{\mu}  \leq c ~ {b}^{\mu+2}.
 \end{align}
% Taking $b=\frac{1}{\theta}$, the assertion holds for the first summation.\\
\end{comment}

 In order to verify the estimation (\ref{main}),  first
 we shall examine the spherical harmonics expansion
 of   $(\cos\theta-1)^N f(\cos\theta)$.  That is, for any  $N\in \NN$, we will find the coefficients $d_l$ in
\begin{align}
(\cos\theta-1)^N f(\cos\theta)= \sum_{l\in \NN_0} d_l~
Z_l(\cos\theta).
\end{align}
(It will follow from our arguments that only the zonal functions
appear  in the
spherical harmonic expansion of function $(\cos\theta-1)^N f$.)\\

Using the  following recursion
 formula   for the Legendre polynomials
 \begin{align}\label{recurtion-formula}
 (2l+1)(x-1)P_l(x)=(l+1)P_{l+1}(x)-(2l+1)P_l(x)+lP_{l-1}(x)\quad \forall x\in [-1,1], ~~l\in \NN_0,
 \end{align}
where we put $P_{-1}\equiv 0$, and with the convention $a_{-1}\equiv
0$, we obtain  the following equalities  for $N=1$:
\begin{align}
(\cos\theta-1)f(\cos\theta)
=& \sum_{l=0} a_l~ (\cos\theta-1)Z_l(\cos\theta)\\
 = & \sum_{l=0}    a_l~ \left\{ \frac{l+1}{2(l+1)+1}Z_{l+1}-\frac{2l+1}{2l+1}Z_l+\frac{l}{2(l-1)+1}Z_{l-1}\right\}(\cos\theta)\\
= &\sum_{l=0} a_l ~ \left\{ \frac{l+1}{2l+3}Z_{l+1}-  Z_l+ \frac{l}{2l-1}~Z_{l-1}\right\}(\cos\theta)\\
=&\sum_{l=1}  a_{l-1}~\frac{l}{2l+1}Z_l(\cos\theta)-\sum_{l=0} a_l~
Z_l(\cos\theta)
+\sum_{l=-1} a_{l+1}~\frac{l+1}{2l+1}Z_l(\cos\theta)\\
=&\sum_{l=0}  a_{l-1}~\frac{l}{2l+1} Z_l(\cos\theta) -\sum_{l=0}
a_l~ Z_l(\cos\theta)
 +\sum_{l=0}a_{l+1}~ \frac{l+1}{2l+1} Z_l(\cos\theta)\\
= &\sum_{l=0}   \left\{
\frac{l}{2l+1}a_{l-1}-a_l+\frac{l+1}{2l+1}a_{l+1}\right\}~
Z_{l}(\cos\theta)\\
= &\sum_{l=0}\frac{l}{2l+1}(a_{l-1}-2a_l+a_{l+1})~Z_l(\cos\theta)  +
\sum_{l=0}\frac{1}{2l+1}(a_{l+1}-a_{l} )~ Z_l(\cos\theta)
\\
%=&\sum_{l=0}\frac{l}{2l+1}(a_{l-1}-2a_l+a_{l+1})~Z_l(\cos\theta)
% + \sum_{l=0}\frac{1}{2l+1}(a_{l+1}-a_{l} )~ Z_l(\cos\theta)\\
%&-(a_1-a_0) Z_0(\cos\theta) - a_0Z_1(\cos\theta)/3
%+a_1Z_0(\cos\theta)\\
%=&\sum_{l=0}\frac{l}{2l+1}(a_{l-1}-2a_l+a_{l+1})~Z_l(\cos\theta)\\
%&+ \sum_{l=0}\frac{1}{2l+1}(a_{l+1}-a_{l} )~ Z_l(\cos\theta)\\
%&+ a_0 Z_0(\cos\theta) -a_0Z_1(\cos\theta)/3\\
%=&\sum_{l=0}\frac{l}{2l+1}(a_{l-1}-2a_l+a_{l+1})~Z_l(\cos\theta)\\
%&+ \sum_{l=0}\frac{1}{2l+1}(a_{l+1}-a_{l} )~ Z_1(\cos\theta)\\
=& \sum_{l=0}a_l^1 ~Z_1(\cos\theta),
\end{align}
where
\begin{align}
a_l^1:= \frac{l}{2l+1}(a_{l-1}-2a_l+a_{l+1})+
\frac{1}{2l+1}(a_{l+1}-a_{l}).
\end{align}

 In fact $a_l^1$ can be written as follows:
\begin{align}
a_l^1:=\left\{R(l)\Delta^+\Delta^- +S(l)\Delta^- \right\}~a_l =
P(l)a_l                     ,
\end{align}
for    the  sequence elements $R(l)$ and  $S(l)$    given as
\begin{align}\label{RSQ}
  R(l)=  \frac{l}{2l+1}, ~~S(l)=\frac{1}{2l+1},
\end{align}
and the operator  $P(l)$    defined as
 $$ P(l):=  R(l)\Delta^+ \Delta^-  +S(l)\Delta^- .$$
 This implies that
  $a_l^1$    comes from  the operation   of  operator $P(l)$ on
   $a_l$; more precisely, say $a:= \{a_l\}$, then

  \begin{align}
  a_l^1 = P(l) a_l =  \left\{R(l)\Delta^+ \Delta^- +S(l)\Delta^-
  \right\}~a_l,
  \end{align}
  and $a_{-1}^l=0$.
  For $N>1$,
analogously, by iteration  one gets
\begin{align}
(\cos\theta-1)^N f(\cos\theta)= \sum_{l=0} a_l^N ~Z_l(\cos\theta),
\end{align}
where $a_l^N=  P(l)^N  (a_l).$ To find a upper  estimate for  the
coefficients $a_l^N$, we shall use the following estimates for the
formulas $R$ and $S$  in (\ref{RSQ}):
 \begin{align}
  R(l) =  \frac{l}{2l+1}= \mathcal{O}(1),
  ~~S(l)=\frac{1}{2l+1}=\mathcal{O}(l^{-1}),
 \end{align}
 and in general:
\begin{align}
(\Delta^-)^{k_1}(\Delta^+)^{k_2}
 R(l) &=  \mathcal{O}(l^{-(k_1+k_2+1)}),\\
  (\Delta^-)^{k_1}(\Delta^+)^{k_2}S(l)&=\mathcal{O}(l^{-(k_1+k_2+1)}).
   %(\Delta^-)^{k_1}(\Delta^+)^{k_2}
% Q(l)&=\mathcal{O}(l^{-(k_1+k_2+2)}).
\end{align}
Using   the preceding  estimates for different exponents $k_i$'s,
the product rules for   difference operators, and   the estimates
in (i) and (ii),
 one obtains the following:
\begin{align}\label{dl}
 a_l^N=  % \{R(l)\Delta^+\Delta^-  +S(l)\Delta^- \}^N a_l &
  \mathcal{O}(l^{ \mu-2N})
  \end{align}
 Now let $N$ be the least integer larger than $\mu/2+1$. The equality (\ref{dl})
implies
\begin{align}
  \mid \sum_l a_l^N ~ Z_l(\cos\theta)\mid   &\leq  c \sum_l  (2l+1) l^{-2N+\mu}\leq C.
\end{align}

  Therefore, since $\mid \cos\theta-1\mid \geq c_1  \theta^2$  for   $0< \theta\leq \pi$ and some constant $c_1$,
   %  for   the second sum in (\ref{spreading})
    we have

\begin{align}\label{second-part}
\left| \sum_l a_l~Z_l(\cos\theta)\right|&=  \frac{1}{\mid
\cos\theta-1\mid^N } \left| \sum_l a_l^N~Z_l(\cos\theta)
\right|\\
& \leq   \frac{c_2}{\mid \cos\theta-1\mid^N }\\
&\leq  \frac{c_3}{  \theta^{2N}},
\end{align}
as desired.

% where $M=-2N+\mu+2$. And,  for the sum over $l< b$ we
%have
%\begin{align}
% \frac{1}{\mid
%\cos\theta-1\mid^N } \left| \sum_{l<b} a_l^N~Z_l(\cos\theta) \right|
 %& \leq \frac{1}{\mid
%\cos\theta-1\mid^N } \sum_{l<b} (2l+1) l^{-2N+\mu} \\
%&  \leq \frac{1}{\theta^{2N} } \sum_{l<b} (2l+1)
%l^{-2N+\mu} \\
%& \leq \frac{1}{\theta^{2N} } \sum_{l<b}  l^2
%l^{-2N+\mu}\\
%&  \leq  c \frac{b^{-2N+\mu+2}}{\theta^{2N} }= c
%\frac{b^M}{\theta^{2N} }
%\end{align}

 %All together,
  % up to a constant independent of $\theta$ and $b$, we have

%\begin{align}
%\mid f(\cos\theta)\mid =  \frac{1}{\mid \cos\theta-1\mid^N } \left|
%\sum_{l=0} a_l^N~Z_l(\cos\theta) \right|\leq  \frac{b^{M}}{\theta
%^{2N} },
%\end{align}
%and hence,  if we set
 %     $b=1/\theta$, we find the estimate of the lemma.\\
 \end{proof}

We are  now ready to state our main result in the next theorem.
 We shall suppose that the sequence $\{c_l\}_{l\in \NN}$ satisfies (i), (ii) of section 3.
%: For $\alpha>2$
%  \begin{itemize}
%\item there exist positive constants $0<C_0\leq C_1< \infty $ such that  for all $l\in \NN$
%$$C_0 l^{-\alpha}\leq  c_l\leq C_1l^{-\alpha},$$
%\item and for any $k_1, k_2\in \NN_0$, $(\Delta^-)^{k_1}{\Delta^+}^{k_2}c_l= \mathcal{O}(l^{-(\alpha+k_1+k_2)})$.
%\end{itemize}
%For example, the $c_l$ could be as in (i), (ii) of section 3.
 Then, we have:
  \begin{theorem}\label{main-theorem} Suppose $\alpha>2$ and $f(s)=s^r f_0(s)$ for a Schwartz function $f_0$ on $\RR^+$ and $r\in \NN$.
Assuming that   $4r+2> \alpha$,    for any   $t>0$ and any two
different points $x, y$ on the sphere
 one has
 \begin{align}\label{corr}
 \left| \frac{ \sum_{l=1} f(t^2\lambda_l)^2 c_l Z_l(x\cdot y)}{ \sum_{l=1}f(t^2\lambda_l)^2 c_l
  Z_l({\bf1})}\right| \leq \frac{Ct^{4r-\alpha+2}}{{d(x, y)}^{2N}},
 \end{align}
 where $N$ is the least integer greater than $2r-\alpha/2 + 1$ and $C$ is independent of $x, y$ and $t$.
 \end{theorem}
 \begin{proof} Choose an interval $I=[a^2,b^2] \subseteq (0,\infty)$ such that
$f^2 \geq c > 0$   on $I$.  Choose $b_1$ with $a < b_1 < b$; then,
if $t$ is sufficiently small, $f(t^2\lambda_l) \geq c$ whenever $a/t
\leq l \leq b_1/t$.  Then, with the assumptions on the coefficients
$c_l)$  we  have the following inequality
 for the  denominator of (\ref{corr}), if $t$ is sufficiently small: Since $Z_l({\bf1}) = (2l+1)$
up to a positive constant, we have
 \begin{align}
  \sum_{l=1}^{\infty} f(t^2\lambda_l)^2 c_l Z_l({\bf 1})&\geq
   \sum_{a/(2t) \leq l\leq b_1/t}f(t^2\lambda_l)^2 c_l Z_l({\bf1})\\
   &\geq
  c~ \sum_{a/t \leq l\leq b_1/t}  \frac{(2l+1) }{l^\alpha  }\\\label{denominator}
  &\geq c~  t^{-1} ~
   t^{-1+\alpha}= c ~t^{\alpha -2} > 0~.
 \end{align}

 To find a upper estimate for the  numerator of (\ref{corr}), for $t>0$,  define

  \begin{align}
 g_t(s)&:=  f(t^2s(s+1))= \left(t^2 s(s+1)\right)^{r} f_0(t^2{s(s+1)}) \quad \text{and}\\
 G_t(s)&:= g_t(s)^2 u(s).
  \end{align}
By induction on $n$, one easily sees that $\partial_s^n
f_0(t^2s(s+1))$ is a finite linear combination of terms of the form
$t^{2i}s^jF(t^2s(s+1))$, where $F \in \mathcal{S}(\RR^+)$, and where
$n=2i-j$.  Since $F \in \mathcal{S}(\RR^+)$, there is a constant $C$
with $|F(t^2s(s+1)| \leq C[t^2s(s+1)]^{-i}$.
%  Recall that $f_0$ in a Schwartz function on $\RR^+$. Therefore, for any $n, m\in \NN_0$ and $t>0$
%there exists a constant $c:=c(n,m)$ such that  for all $s>0$
%   \begin{align}
%    |  \partial_s^n f_0(t^2s(s+1))|& \leq c ~ t^{2n} s^n( t^2s(s+1))^{-m}\\
%    &\leq t^{2(n-m)}\min\{ s^{n-2m}, s^{n-m}\}
%  \end{align}
 Hence, for  $ s\geq 1$
  \begin{align}
      | \partial_s^n f_0(t^2s(s+1))|\leq c ~  s^{-n}.
   \end{align}
%This  implies that for any $n\in \NN_0$ there exists a constant $c:=c_n$ that
%      \begin{align}
%    |  \partial_s^n f_0(t^2s(s+1))^2|\leq c ~  s^{-n} \quad ~\forall ~ s\geq 1
%  \end{align}
%holds uniformly.
 % Using the preceding estimation for $f_0$,
 Thus, for    $i \in \NN_0$   the inequality
  \begin{align}\label{est-on-g}
  | \partial^i g_t(s)^2|\leq  t^{4r}s^{4r-i} \quad \forall ~ s\geq 1
  \end{align}
  holds uniformly  up to a constant which is   independent of $t$.
  The  estimation (\ref{est-on-g}) and the assumptions on  $u$,   that for any $j\in \NN_0$
the inequality $|\partial^j u(s)|\leq c~s^{-\alpha-j}$  holds
uniformly up to a constant $c$,  imply that  for $t>0$ and  $k\in
\NN_0$   the estimation
    \begin{align}
 | \partial^k G_t (s)|\leq   t^{4r}s^{4r-k-\alpha}
  \end{align}
  holds  uniformly up to a constant, independent of the parameter $t$.
%Therefore for   $t>0$
     %for $i, k\geq 1$

By Lemma \ref{tool}, for $\mu=4r-\alpha$ we have
% \begin{align}\label{ohne-t}
%  \left|\sum_{l=1} G(l)  Z_l(x \cdot y)\right|&=
%  \left|\sum_{l=1} f(\lambda_l)^2 c_l Z_l(x\cdot y)\right|\\
% & \leq \frac{c}{d(x , y) ^{4r-\alpha+2}}.
% \end{align}
%and hence
   \begin{align}\label{mit-t}
    \left|\sum_{l=1} G_t(l)  Z_l(x \cdot y)\right|&=
 \left|\sum_{l=1} f(t^2\lambda_l)^2 c_l Z_l(x\cdot y)\right|\\
  &\leq
  ~c~\frac{t^{4r}}{\left(d(x, y)\right)^{2N}}.
%& = ~ c~ \frac{ t^{\alpha-2}}{\left(\frac{d(x,
%y)}{t}\right)^{2N}}.
 \end{align}
Using the estimate (\ref{denominator}) for the denominator, we
therefore have that

\begin{align}
 \left| \frac{ \sum_{l=1} f(t^2\lambda_l)^2 c_l Z_l(x\cdot y)}
 { \sum_{l=1}f(t^2\lambda_l)^2 c_l Z_l({\bf1})}\right| \leq \frac{t^{4r-\alpha+2}}{d(x, y)^{2N}},
 \end{align}
as desired.
 \end{proof}

 \begin{remark}
 Note that the
  result of Theorem \ref{main-theorem} is  true for $r=1$ if $\alpha< 6$.
 \end{remark}

\end{document}